\documentclass[12pt]{article}
\usepackage{amsmath,latexsym}
\usepackage{amsthm}
\usepackage{amssymb}
\usepackage{amsfonts}
\usepackage{graphicx}

\title{Solutions of the problem of Erd\"os-Sierpi\'nski: $\sigma(n)=\sigma(n+1)$ }
\author{Lourdes Benito\\
\small{lourdes.benito@estudiante.uam.es}}
\date{\centering }

\begin{document}

\maketitle

\newtheorem{theorem}{Theorem}

\newtheorem{problem}{Problem}
\begin{abstract}
For $n\leq 1.5 \cdot 10^{10}$, we have found a total number of 1268
solutions
 to the Erd\"os-Sierpi\'nski problem  finding positive
integer solutions of $\sigma(n)=\sigma(n+1)$, where $\sigma(n)$ is
the sum of the positive divisors of n. On the basis of that set of
solutions the following empirical properties are enunciated: first,
all the $\sigma(n)$, $n$ being a solution, are divisible by $6$;
second, the repetition of solutions leads to the formulation of a
new problem: \emph{ Find the natural numbers $n$ such that
$\sigma(n)=\sigma(n+1)=\sigma(n+k)=\sigma(n+k+1)$ for some positive
integer $k$}. A third  empirical property concerns the  asymptotic
behavior of the function of $n$ that gives the number of solutions
for $m$ less or equal to $n$, which we find to be as $n^{1/3}$.
Finally some theorems related to the Erd\"os-Sierpi\'nski problem
are enunciated and proved.
\end{abstract}

\section{Introduction}

 The
problem of determining the natural numbers $n$ that verify
\begin{equation}\label{sigma}
\sigma(n)=\sigma(n+1)
 \end{equation}
 Where  $\sigma(n)$ is the sum of the positive divisors of $n$,
 is known as the
Erd\"os-Sierpi\'nski problem.

Erd\"os \cite{Erd} had conjectured that there are infinitely many
solutions to the (\ref{sigma}) but he offered no proof. Some time
later, Sierpi\'nski  \cite{Sier} posited  the same question.

The first natural numbers that satisfy  \eqref{sigma} are:
$$n=14, \  206,\  957, \ 1334, \ 1364,\  1634,\  2685\ldots$$

This sequence is identified as the  A002961 sequence by the Sloane's
On-Line Encyclopedia of Integer Sequences.

Guy \cite{Guy} , B13,   communicates that   Jud McCranie  calculated
solutions for  the problem of Erd\"os-Sierpi\'nski up to
$n=4.25\cdot  10^9 $. A total of $832$ solutions were found. We have
investigated the solutions up to  $n=1.5 \cdot 10^{10}$ (see \ref{calculo}), and found a total of
 1268 solutions.

Guy and Shanks \cite{GS}  observe  that the solutions for $n=14$,
$n=206$ and $n=19358$
 are given by:
\begin{equation}\label{gsfor1}
n=2p, \quad n+1=3^m q,
\end{equation}
 where
$$q=3^{m+1} -4, \quad p=(3^m q  -1)/2$$
are both prime, and $m$ equals $1$, $2$ or $4$.

They also indicate that  the solutions for  $n=18873$, $n=174717$
and $n=5559060136088313$ are given by
\begin{equation}\label{gsfor2}
 n=3^m q, \quad n+1= 2p,
\end{equation}
with the primes
$$q=3^{m+1} -10 , \quad p=(3^m q +1)/2$$
for $m=4$, $5$ and $16$.

With the help of the two functions {\tt guyshanks1(1,1000)} and {\tt
guyshanks2(1,1000)}, written in Pari-GP language \cite{pari}  (see \ref{apend_gsfor}), it is  verified that  no other solutions of the
forms \eqref{gsfor1} and \eqref{gsfor2} exist for $m\leq 1000$.

\section{Some Empirical Results}

On the basis of the set of solutions of the Erd\"os-Sierpi\'nski
problem for  $n \le 1.5 \cdot 10^{10}$, the following empirical
results can be enunciated:

For $n$ a solution to the Erd\"os-Sierpi\'nski problem

\begin{itemize}
\item[i)]\label{i)}  $\sigma(n)$ is divisible by $6$.
\item[ii)]$\sigma(n)$  is divisible by $4$ with the exception of:
$$\sigma(18873)=  28314= 2   \cdot  3^2 \cdot   11^2 \cdot   13$$
\item[iii)] $\sigma(n)$ is divisible by $8$ except the mentioned
case ii) at the following ones:
$$\sigma( 4364 )= 7644 = 2^2\cdot      3 \cdot    7^2 \cdot 13 $$
$$\sigma(  14841 )=  22932  =  2^2 \cdot 3^2 \cdot   7^2 \cdot   13 $$
$$\sigma(  3582224 )=  6976860 =  2^2\cdot   3 \cdot    5\cdot    11^2\cdot        31^2$$
$$\sigma(195694137  )=  293616180 =  2^2\cdot     3^2\cdot     5\cdot    11^2\cdot    13\cdot    17\cdot    61$$
$$\sigma(597311577)=  896138100=  2^2\cdot     3^3\cdot     5^2\cdot 11^2\cdot 13\cdot   211 $$
\item[iv)]The only $\sigma(n)$  of the form $2^a \cdot 3^b$ with
$a$ and $b$ natural numbers are:
$$\sigma(14)=24=2^3\cdot3$$
$$\sigma(147454)= 221184=2^{13}\cdot3^2$$
\end{itemize}

\subsection{A new problem}\label{s_nuevoproblema}

If the solutions to the Erd\"os-Sierpi\'nski problem are carefully
examined, it will be observed that certain values are repeated.
Therefore, a related problem can be enunciated:

\begin{problem}\label{problemanuevo1}
 Find the  natural numbers $n$ that  verify
\label{problemanuevo}
$$\sigma(n)=\sigma(n+1)=\sigma(n+k)=\sigma(n+k+1)$$
for some positive integer $k$.
\end{problem}

 For $n\le 1.5 \cdot 10^{10}$,  22 solutions  to this  problem have been found (see \ref{solrep}).
The smallest  solutions is:
$\sigma(79826)=\sigma(79827)=\sigma(79833)=\sigma(79834)=120960$.

$\sigma(n)=4049740800$ has three different solutions. In this case,
the problem can enunciated in a different way:

\begin{problem}\label{problemanuevo2}
Find the  natural numbers $n$ that  verify
$$\sigma(n)=\sigma(n+1)=\sigma(n+k)=\sigma(n+k+1)=\sigma(n+k')=\sigma(n+k'+1)$$
for some pair of positive integer numbers  $k$ and $k'$.
\end{problem}

\subsection{Estimation of the number of solutions to the Erd\"os-Sierpi\'nski problem}\label{s_numerosoluciones}

\begin{figure}[h]
\centering
\includegraphics[width=8cm]{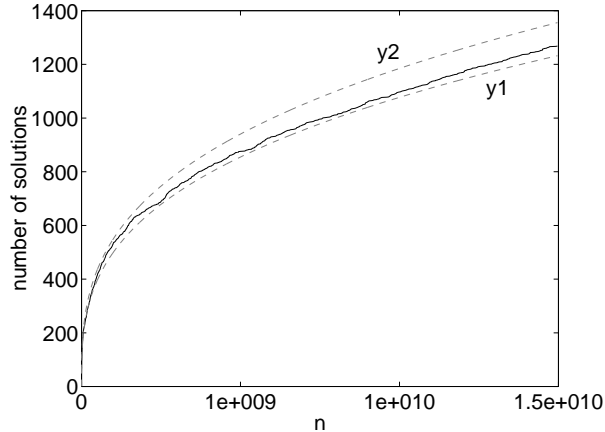}\\
\caption{Number of solutions for (\ref{sigma}) less or equal to $n$ versus
$n$. Dashed lines correspond to $y_1=0.5 \cdot n^\frac{1}{3}$ and
$y_2=0.55 \cdot n^{\frac{1}{3}}$. The maximum value of $n$ for which
the number of solutions is over $y_1$ is $n=258083942$, and the
maximum value of $n$ for which the number of solutions is below
$y_2$ is $n= 2305557$. \label{fig1}}
\end{figure}

\begin{figure}[h]
\centering
\begin{tabular}{cc}
\includegraphics[width=6.5cm]{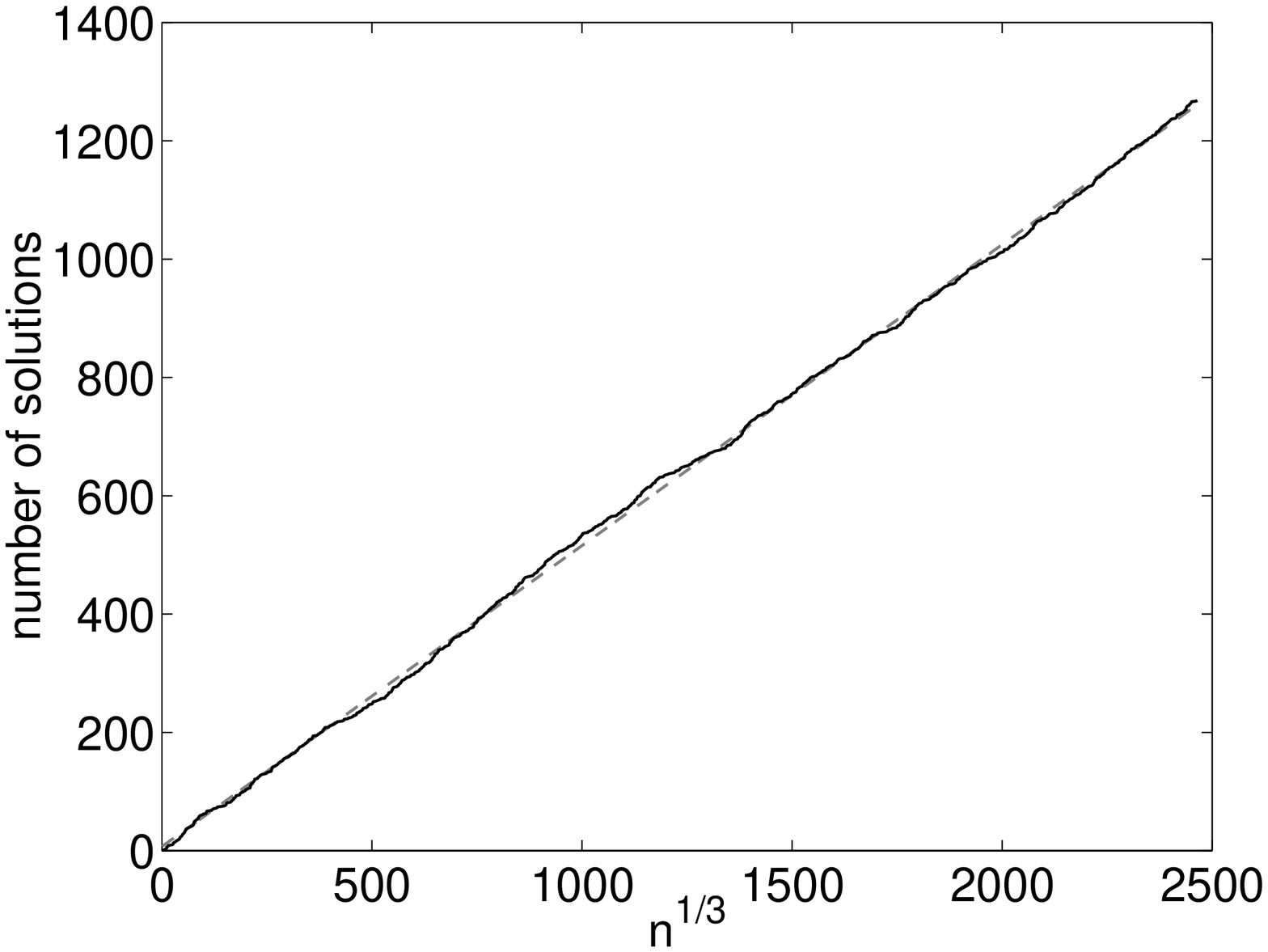}&\includegraphics[width=6.5cm]{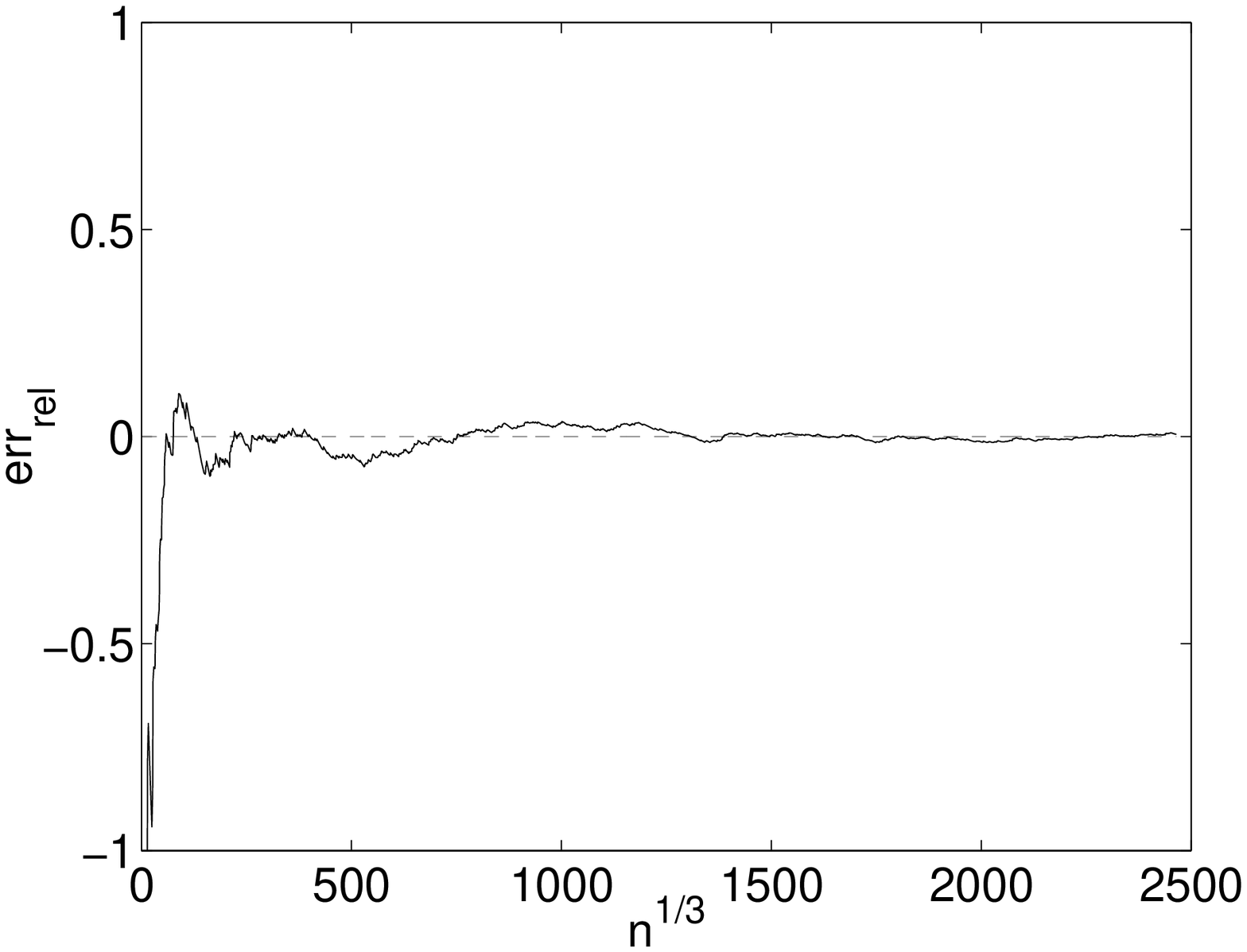}\\
\end{tabular}
\caption{On  the left, both the number of solutions less or equal to
$n$ (continuous line) and the linear adjustment (dashed line) are
plotted versus $n^{1/3}$. On the right, it is represented the
relative error ($\epsilon$) corresponding to the adjustment.
$\epsilon=\frac{y-y_{adj}}{y}$, where $y$ is the number of solutions
less or equal to $n$ and $y_{adj}$ is the
 number of solutions that predicts the adjustment.
From  $n=792855$ $(n^{\frac{1}{3}}= 92.5546)$ the relative error is
smaller than  $10\%$ $(|\epsilon|<0.1)$.\label{fig2yfig3}}
\end{figure}

In order to give an estimation of the number of solutions, different
fit functions were tried, and the best fit obtained was  in the form
$\sim n^{\frac{1}{3}}$.

Then, it is found that, for $258083942<n<1.5\cdot 10^{10}$ the
number of solutions lies between   $y_1$ and $y_2$

$$
\left\{
\begin{array}{lcr}
y_1&=&0.50 \cdot n^{\frac{1}{3}}\\
 y_2&=&0.55 \cdot n^{\frac{1}{3}}
\end{array}
\right.
$$

If a linear adjustment based in $n^{\frac{1}{3}}$ is made, we get:

$$y_{adj}= 0.5088 \cdot n^{\frac{1}{3}} + 6.9183$$
\noindent where $y_{adj}$ is the number of solutions less or equal
to $n$ that predicts the adjustment (see figure \ref{fig2yfig3}).

For example, an estimation of the number of solutions up to
$n=10^{15}$ to the  Erd\"os-Sierpi\'nski problem is: $ 0.5088
\cdot 10^{15\cdot \frac{1}{3}}+ 6.9183 \approx 5 \cdot 10^{4} $
solutions.

\section{Some theorems}

\begin{theorem}\label{teorema1}
$\sigma(n)$ is odd  if and only if $n$ is a square, or $n$ is the
double of a square.\end{theorem}

\begin{proof}
\

\begin{itemize}
\item If  $p$ is prime and $p\not= 2$ $\Rightarrow$
$\sigma(p^{\alpha})=1+p+p^2+\ldots+p^{\alpha} \equiv \alpha+1 \mod
2$, then for  $\sigma(p^{\alpha})$ to be odd, $\alpha$ must be
even. \item If $p=2$ $\Rightarrow$ $\sigma(2^{\alpha})= 1+ 2 +
\ldots + 2^{\alpha}$ is always odd.
\end{itemize}

Then for $\sigma(n)$  to be odd, in the decomposition in prime
factors of $n$, the power of the primes $p>2$ must be even.
Therefore, $n$ is a square number or a $2\cdot$square number.
\end{proof}

\begin{theorem}

For every  natural number $n$ that verifies $\sigma(n)=\sigma(n+1)$,
 $\sigma(n)$ is multiple of $2$ or $3$.
\end{theorem}

\begin{proof}

If $\sigma(n)$ is odd, from theorem \ref{teorema1} and, since two
consecutive positive integer numbers cannot be both square,
 either $n$ must be square and $n+1$
is the double of an square or $n$ is the double of a  square and
$n+1$ is an square.

Then if   $\gcd (x,2)=1$, $\sigma(2^{2\alpha+1}\cdot
x^2)=\sigma(2^{\alpha+1}) \cdot \sigma(x^2)= (1+2+2^2+\ldots  +
2^{2\alpha+1}) \cdot \sigma(x^2) = (1+2+2^2(1+2)+\ldots +
2^{2\alpha}(1+2))\cdot \sigma(x^2)=3\cdot (1+2^2 +2^4 + \ldots +
2^{2\alpha})\cdot \sigma(x^2)$.
\end{proof}

\begin{theorem}\label{teoremapell}
If  $n$ is a positive integer such that $\sigma(n)=\sigma(n+1)$ is
odd, then:

$$n=2 y^2,\  n+1=x^2 \hbox{ with } x^2-2y^2=1$$
or
$$n=x^2,\ \  n+1=2 y^2 \hbox{ whit }  x^2-2y^2=-1, $$
where $x$ and $y$ are integer numbers.
\end{theorem}

\begin{proof}
From  to theorem \ref{teorema1} if $\sigma(n)$ is odd, then $n$ must
be square or $2\cdot$square. If $n$ is odd, then $n=x^2$, and
$n+1=2y^2$  because two consecutive positive integers   cannot be
both square numbers. Therefore  $x^2-2y^2=-1$.

If $n$ is even then  $n=2y^2$, and  $n+1=x^2$. Therefore
$x^2-2y^2=1$.
\end{proof}

Using a Pari-Gp \cite{pari} function (see \ref{sigimp})  based on
the ideas of \cite{brocot} we have verified that there is not any
solution to the problem of Erd\"os-Sierpi\'nski in which
$\sigma(n)$ is odd for $n<10^{50}$.

\section{Open questions}

\begin{enumerate}
\item The  numerical evidences  described in section
\ref{s_numerosoluciones},  agree with the conjecture of Erd\"os
about  the existence of infinite $n$   that satisfy
$\sigma(n)=\sigma(n+1)$.
 These evidences allow us  to conjecture that the number of solutions less or equal  than $n$ is of the order of $\displaystyle  n^{\frac{1}{3}}$.

\item Apart from those mentioned, do any other solutions to the
Erd\"os-Sierpi\'nski problem  of forms \eqref{gsfor1} and
\eqref{gsfor2} exist?

\item Does  any $n$ for which
$\sigma(n)=\sigma(n+1)\not=\dot{2}$ exist?

\item Does  any $n$ for which
$\sigma(n)=\sigma(n+1)\not=\dot{3}$ exist?

\item Apart from those mentioned, do any other solutions to the
Erd\"os-Sierpi\'nski problem in the  form
$\sigma(n)=\sigma(n+1)=2^a\cdot3^b$ exist?

\item Do  infinite solutions to problem
\ref{problemanuevo1} exist?

\item Apart from those mentioned, do any other solutions to the problem
\ref{problemanuevo2} exist?

\end{enumerate}

\appendix

\section{Numerical code used in section \ref{gsfor1}}\label{apend_gsfor}

Below the reader will find the code of the  Pari-Gp function that
verifies that there are not more solutions of form (\ref{gsfor1})
for values of $m$  less than 1000.




\begin{verbatim}
q(m)=3^(m+1)-4

p(m)=(3^(2*m+1)-4*3^m-1)/2

guyshanks1(k1,k,i)=
     {
      for(i=k1,k,
        if(isprime(q(i)),
             if(isprime(p(i)), print(i," ",2*p(i))," ",3^i*q(i))
           ,)
          )
      }
\end{verbatim}

\newpage

And the code for the  Pari-Gp function that verifies that there are
not more solutions of form(\ref{gsfor2}) for  values of $m$ less
than 1000.


\begin{verbatim}
q(m)=3^(m+1)-10

p(m)=(3^(2*m+1)-10*3^m+1)/2

guyshanks2(k1,k,i)=
   {
    for(i=k1,k,
      if(isprime(q(i)),
         if(isprime(p(i)),print(i," ",3^i*q(i)," ",2*p(i)),)
         ,)
        )
   }
\end{verbatim}

\section{Numerical code  for theorem 3}\label{sigimp}

\begin{verbatim}
noteven(k,z,t,a)=
    { z=1;
       while(denominator(z)<k,z=2+1/z;t=1+1/z;
           a=sigma(denominator(t)^2);
              if(a==sigma(numerator(t)^2\2),
                print(denominator(t)^2," ",numerator(t)^2\2," ",a),)
             )
     }
\end{verbatim}

\section{Solutions for the problems 1 and 2}\label{solrep}

\footnotesize

\begin{tabular}{|c|c|c|c|c|c|}
\hline
$\sigma(n)$&$n$&index $n$&$n+k$&index $n+k$& $k$\\
\hline\hline
120960 &79826 &21 &79833 &22 &7 \\
\hline
17418240 &9669915 &112 &10195305 &116 &525390 \\
\hline
27095040 &13577355 &129 &17714486 &135 &4137131 \\
\hline
220631040 &105967784 &236 &109606496 &240 &3638712 \\
\hline
309657600 &157073930 &266 &186027314 &286 &28953384 \\
\hline
339655680 &177961544 &279 &208193654 &296 &30232110 \\
\hline
840153600 &450427130 &399 &552488553 &430 &102061423 \\
\hline
1132185600 &610698123 &448 &628487810 &453 &17789687 \\
\hline
2868203520 &1593132723 &619 &1717235637 &635 &124102914 \\
\hline
3096576000 &1974224366 &652 &2038559845 &659 &64335479 \\
\hline
3690086400 &2004024915 &654 &2008900875 &655 &4875960 \\
\hline
4049740800 &2155986584 &667 &2482257261 &686 &326270677 \\
\hline
4049740800 &2155986584 &667 &2663065825 &712 &507079241 \\
\hline
4049740800 &2482257261 &686 &2663065825 &712 &180808564 \\
\hline
4790016000 &2750303997 &725 &3071446581 &749 &321142584 \\
\hline
6991246080 &3639165255 &795 &4642434141 &857 &1003268886 \\
\hline
10879989120 &5479093304 &893 &6209887131 &938 &730793827 \\
\hline
11544519168 &7038224205 &979 &7566526345 &998 &528302140 \\
\hline
11805393600 &6421134644 &951 &7559982333 &997 &1138847689 \\
\hline
12773376000 &8053053458 &1014 &8334008133 &1024 &280954675 \\
\hline
15926803200 &9660726662 &1079 &10458788337 &1112 &798061675 \\
\hline
21250252800 &11017504952 &1137 &13099733673 &1208 &2082228721 \\
\hline
\end{tabular}

\normalsize

\section{Solutions for the problem of \hbox{Erd\"os-Sierpi\'nski} up to $n=1.5\cdot10^{10}$}\label{calculo}

\scriptsize
%

\end{document}